\RequirePackage{fix-cm}
\documentclass[smallextended,envcountsect]{svjour3}       % onecolumn (second format)
\smartqed  % flush right qed marks, e.g. at end of proof
\usepackage{microtype} % micro layout optimization
\usepackage{graphicx}
\usepackage{mathtools}
\usepackage{stmaryrd}
\usepackage{mathrsfs}
\usepackage{booktabs}
\usepackage{multirow}
\usepackage{amssymb}
\usepackage[title]{appendix}
\usepackage{xcolor}
\usepackage{textcomp}
\usepackage{manyfoot}
\usepackage{algorithm}
\usepackage{algorithmicx}
\usepackage{algpseudocode}
\usepackage{listings}
\usepackage{hyperref}
\usepackage{cleveref}
\newcommand{\R}{\mathbb{R}}
\renewcommand{\d}{\mathrm{d}}
\newcommand{\norm}[1]{\left\lVert#1\right\rVert}
\newcommand{\abs}[1]{\left\lvert#1\right\rvert}
\newcommand{\jump}[1]{\left\llbracket#1\right\rrbracket}
\newcommand{\avg}[1]{\left\{#1\right\}}
\DeclareMathOperator{\Grad}{grad}
\DeclareMathOperator{\Div}{div}
\DeclareMathOperator{\diag}{diag}

\numberwithin{equation}{section}
\journalname{Journal of Scientific Computing}
\begin{document}

\title{Unconditionally $L^2$-stable explicit--implicit--null filtered Runge--Kutta methods for general nonlinear degenerate parabolic equations%\thanks{Grants or other notes
}

\titlerunning{Unconditionally $L^2$-stable EIN filtered RK}        % if too long for running head

\author{Yue Wu \and Chi-Wang Shu}

\authorrunning{Y. Wu \and C.-W. Shu} % if too long for running head

\institute{Y. Wu \at
    Division of Applied Mathematics, Brown University, Providence, RI 02912, USA. \\
    \email{yue\_wu3@brown.edu}
    \and
    C.-W. Shu \at
    Division of Applied Mathematics, Brown University, Providence, RI 02912, USA. \\
    \email{chi-wang\_shu@brown.edu}
}

\date{Received: date / Accepted: date}

\maketitle

\begin{abstract}
We develop unconditionally $L^2$-stable explicit--implicit--null (EIN) filtered Runge--Kutta (RK) methods for nonlinear degenerate parabolic equations with space-, time-, state-, and gradient-dependent diffusion coefficients. The methods are linearly implicit in the sense that they evaluate nonlinear terms explicitly and require only linear solves. Unlike classical EIN methods that are direct applications of implicit--explicit (IMEX) RKs, our construction is built upon explicit RKs. The key ingredient is two-sided symmetric filtering applied at each RK stage: a high-order filter is applied both before and after evaluating the original nonlinear residual. It leads to a dissipation estimate independent of the time step and retains the underlying RK order formally. Assuming time-dependent splitting operators with stagewise pointwise coefficient bounds, we rigorously prove non-increasing $L^2$ norm for arbitrary time steps. A quadratic matrix criterion covers explicit RK-based methods with positive weights, while strong stability preserving (SSP) compositions yield exact energy identities and weighted dissipation estimates. We present SSP methods through fourth order and a method based on the classical RK4. The framework accommodates spatial discretizations with discrete gradient structures including continuous Galerkin, Bassi--Rebay 1, the method of Brezzi et al., and local discontinuous Galerkin. Compatible summation-by-parts discretizations with simultaneous approximation terms are also covered in their discrete mass norm.
\keywords{explicit--implicit--null time stepping
\and Runge--Kutta, 
\and strong stability preserving, 
\and degenerate parabolic equations, 
\and nonlinear diffusion, 
\and linearly implicit time discretization}
\subclass{65M12 \and 65L06 \and 65L04 \and 35K65}
\end{abstract}

\section{Introduction}\label{sec:introduction}

\subsection{Explicit--implicit--null time stepping}

Explicit time stepping for diffusion is subject to a parabolic time step restriction \(\Delta t = \mathcal{O}(h^2)\). A fully implicit treatment avoids this restriction at the cost of complicated nonlinear iterations. The idea of explicit--implicit--null (EIN) dates back to the alternating-direction Galerkin method of Douglas and Dupont \cite{10.1016/B978-0-12-358502-8.50009-8}, and the framework is named by Duchemin and Eggers \cite{10.1016/j.jcp.2014.01.013}. It combines both advantages of explicit and implicit time discretizations by adding and subtracting a linear damping operator to the right-hand side of the PDE, treating the added term implicitly and the remainder explicitly. Hence it can be called \textit{linearly implicit}. Explicit evaluation of the nonlinear terms avoids Newton--Raphson or Picard iterations, while implicit treatment of a suitable linear operator expects unconditional stability. EIN fits within the IMEX Runge--Kutta (RK) framework \cite{10.1137/0732037,10.1016/S0168-92749700056-1,10.1007/s10915-004-4636-4} and easily extends to high-order time stepping. For PDEs involving nonlinear diffusion, EIN time stepping can be coupled with a wide range of spatial discretizations including finite difference \cite{10.1016/j.jcp.2014.01.013,10.1016/j.jcp.2022.111314,10.1016/j.jcp.2024.113626,10.1016/j.compgeo.2026.108315}, local discontinuous Galerkin (LDG) \cite{10.1007/s11425-018-9524-x,10.1016/j.jcp.2022.111314,10.1007/s42967-023-00359-6}, spectral collocation \cite{10.4208/eajam.2022-271.090123}, direct DG (DDG) \cite{10.4208/nmtma.OA-2024-0076}, and also supports a natural extension to variable splitting coefficients \cite{10.1007/s42967-023-00359-6}. For Allen--Cahn-type equations, EIN can also use point-wise damping as the splitting and employ efficient ETDRK or IFRK \cite{10.1016/j.cma.2022.114817,10.1016/j.apnum.2023.03.002}. To the authors' best knowledge, EIN methods for nonlinear diffusion problems have the following two limitations. 
\begin{enumerate}
	\item The splitting operator is key to the stability of the EIN methods, but there has been few rigorous \(L^2\) stability analyses on fully discrete EIN methods for \textit{general} nonlinear degenerate diffusion. Most of the stability analyses are von Neumann analyses applied to constant-coefficient linear cases, including \cite{10.1016/j.jcp.2014.01.013,10.1007/s11425-018-9524-x,10.1016/j.jcp.2022.111314,10.4208/eajam.2022-271.090123,10.4208/nmtma.OA-2024-0076,10.1016/j.compgeo.2026.108315}. Semi-discrete \(L^2\) stability analyses include \cite{10.1007/s42967-023-00359-6}. The fully discrete one in \cite{10.1016/j.cnsns.2024.108094} requires separable diffusion coefficients and small time steps and results in a Gr\"onwall-type \(L^2\) estimate, and our recent work in \cite{10.48550/arXiv.2610.34539} is limited to a special case which admits an energy potential and does not imply \(L^2\) stability in the degenerate case. 
	\item All of the relevant EIN methods are direct applications of additive IMEX-RK methods \cite{10.1137/0732037,10.1016/S0168-92749700056-1,10.1007/s10915-004-4636-4} that require two Butcher tableaux. Construction of EIN time stepping using general explicit RKs remains an open problem. 
\end{enumerate}

For nonlinear degenerate parabolic equations with a gradient-flow structure, Yan et al. \cite{10.1016/j.cam.2023.115674} study entropy-dissipative implicit LDG schemes with a KKT limiter for positivity and mass conservation. They prove entropy dissipation for implicit Euler without a positivity constraint and analyze the positivity-constrained first-order scheme under a time-step restriction.

\subsection{The model problem and \texorpdfstring{$L^2$}{L2} stability}

Consider the nonlinear degenerate parabolic problem
\begin{subequations}\label{eq:model}
	\begin{alignat}{2}
		u_t - \Div \left(a(t, \mathbf{x}, u, \Grad u) \Grad u \right) & = 0, & \quad & (t, \mathbf{x}) \in (0,T] \times \Omega, \label{eq:original-PDE} \\
		u(0, \mathbf{x}) & = u_0(\mathbf{x}), & & \mathbf{x} \in \Omega, \label{eq:original-IC}\\
		u(t, \mathbf{x}) & = 0, & & (t, \mathbf{x}) \in (0,T] \times \partial\Omega. \label{eq:original-Dirichlet-BC}
	\end{alignat}
\end{subequations}
Here $\Omega \subset \R^d$ is a bounded Lipschitz domain and $a: \R \times \R^d \times \R \times \R^d \mapsto [0, +\infty)$ may vanish. The Dirichlet BC \eqref{eq:original-Dirichlet-BC} can also be replaced with periodic or homogeneous Neumann boundary conditions. Testing with $u$ gives, for a sufficiently regular solution,
\begin{equation}\label{eq:semi-discrete-stability}
	\frac{1}{2}\frac{\d}{\d t}\norm{u}_{L^2(\Omega)}^2
	= -\int_\Omega a(t,\mathbf{x},u,\Grad u)\abs{\Grad u}^2\,\d\mathbf{x}
	\leq 0,
\end{equation}
and consequently,
\begin{equation}\label{eq:space-time-stability}
	\frac{1}{2} \norm{u(T)}_{L^2(\Omega)}^2 + \int_0^T \int_\Omega a(t,\mathbf{x},u,\Grad u)\abs{\Grad u}^2 \, \d\mathbf{x} \, \d t 
	= \frac{1}{2} \norm{u_0}_{L^2(\Omega)}^2.
\end{equation}
We seek fully discrete counterparts of \eqref{eq:semi-discrete-stability} that give one-step strong \(L^2\) stability \(\norm{u_h^{n+1}}_{L^2(\Omega)} \leq \norm{u_h^{n}}_{L^2(\Omega)}\) or \eqref{eq:space-time-stability} that control $\norm{u_h^{n}}_{L^2(\Omega)}$ and the accumulated dissipation by $\norm{u_h^0}_{L^2(\Omega)}$ uniformly in the time step \(n\).

The EIN splitting of \eqref{eq:original-PDE} is 
\begin{equation}\label{eq:split-intro}
	u_t = L(t,u) = L^E(t,u) + L^I(t,u)
\end{equation}
where 
\begin{subequations}
	\begin{align}
		L(t,v) & \coloneqq \Div(a(t,\mathbf{x},v,\Grad v) \Grad v), \\
		L^E(t,v) & \coloneqq S(t)v - \mathcal{A}(t,v)v, \\
		L^I(t,v) & \coloneqq -S(t)v, \\
		\mathcal{A}(t,v) & \coloneqq -\Div(a(t,\mathbf{x},v,\Grad v) \Grad\cdot), \\
		S(t) & \coloneqq -\Div(s(t,\mathbf{x}) \Grad\cdot)
	\end{align}
\end{subequations}
Here $\mathcal A(t,v)$ is the frozen self-adjoint semi-positive nonlinear diffusion operator evaluated at some \(v\), and $S(t)$ is the self-adjoint semi-positive splitting operator with a coefficient $s: \R \times \R^d \mapsto [0, +\infty)$. This formulation \eqref{eq:split-intro} belongs to the variable-coefficient EIN splitting \cite{10.1007/s42967-023-00359-6}, where the splitting coefficient may vary in both space and time and is selected separately at each RK evaluation stage.

\subsection{Our approach and contributions}

For general diffusion coefficients, $\mathcal A(t,v)$ and $S(t)$ do not necessarily commute. Even when a sufficiently large splitting coefficient gives $0 \preceq \rho\mathcal A(t,v) \preceq S(t)$ for some $\rho > 0$, this comparison alone does not guarantee $L^2$ contraction of a one-sided resolvent update; see the counterexample in \cref{sec:motivation}. Motivated by the two-sided resolvent stability \cref{prop:two-sided-ein1}, our construction applies the same high-order filter before and after evaluating the nonlinear diffusion residual. The resulting symmetric product involving the frozen diffusion operator gives a stagewise bound relating dissipation to the squared residual, uniformly in the time step and even the underlying spatial discretization.

This bound yields EIN filtered methods with a non-increasing $L^2$ norm for arbitrary time steps, constructed from a single explicit RK formula. To the best of our knowledge, this construction differs from all previous EIN time discretizations, including those in our recent preprint \cite{10.48550/arXiv.2610.34539}: instead of directly applying an additive IMEX-RK method with two Butcher tableaux, we build high-order methods by combining an explicit RK formula with symmetric filtering. A quadratic matrix criterion covers any explicit RK method with positive weights (\cref{cor:positive-rk}), while the convex decomposition of SSP methods into forward Euler branches \cite{10.1137/S003614450036757X} gives exact energy identities and weighted physical and mixing dissipation estimates (\cref{thm:ssp-exact-general}). The matrix criterion complements existing RK energy analyses of linear seminegative systems \cite{10.4310/AMSA.2017.v2.n2.a3,10.1137/18M122892X,10.1137/22M1472218} including the possible failure of one-step norm contraction for classical RK4. Here the uniform filtered-stage bound permits nonlinear frozen diffusion operators that differ between stages. The methods evaluate nonlinear terms explicitly and require only symmetric positive definite linear solves. They formally retain the temporal order of the underlying RK formula.

We give explicit SSP-EIN1--3 formulas and obtain their dissipation estimates from the general SSP theorem, and also construct methods based on SSPRK(5,4) and classical RK4. For SSP-EIN1--3, a sufficient stagewise splitting condition is
\begin{equation}\label{eq:stage-bound-intro}
	0 \leq C_p^{\mathrm{EIN}} a(t^{n,(i)}, \mathbf{x}, w_h^{n,(i)}(\mathbf{x}), G_h w_h^{n,(i)}(\mathbf{x})) 
	\leq s^{n,(i)}(\mathbf{x}),
\end{equation}
where $w_h^{n,(i)}$ is the filtered evaluation stage and $G_h w_h^{n,(i)}$ its discrete gradient. The sufficient factors are
\begin{equation}\label{eq:thresholds-intro}
	C_1^{\mathrm{EIN}} = \frac{1}{4}, \quad
	C_2^{\mathrm{EIN}} \approx 0.619684, \quad
	C_3^{\mathrm{EIN}} \approx 1.010104, 
\end{equation}
with the general bound $C_p^{\mathrm{EIN}} < \frac{p^2}{4}$ for \(p > 1\) established in \cref{lem:sharp-resolvent}. The stagewise assumptions allow degenerate diffusion and noncommuting operators, including splitting operators that vary between stages. General positive-weight RK formulas require the additional RK factor in \eqref{eq:positive-rk-bound}.

The spatial framework covers continuous Galerkin (CG) \cite{10.1016/S0168-2024_08_70184-X}, Bassi--Rebay 1 (BR1) \cite{10.1006/jcph.1996.5572}, the method of Brezzi et al. \cite{BrezziEtAl1999Diffusion}, and local discontinuous Galerkin (LDG) \cite{10.1137/S0036142997316712}. Their discrete gradients and adjoint divergences provide a common basis for the stability proof \cite{10.1016/j.cam.2015.10.024}, with DG lifting forms described in \cite{10.1137/S0036142901384162}. The coefficient-weighted realizations use a fixed positive definite mass product and symmetric semi-positive frozen diffusion forms, and allow arbitrary polynomial degree. Compatible summation-by-parts (SBP) discretizations with simultaneous approximation terms (SATs) \cite{10.1016/j.jcp.2021.110634} also fit this framework when the complete stabilized operator admits the required gradient--adjoint factorization.

Our recent preprint \cite{10.48550/arXiv.2610.34539} studies state-dependent diffusion through a Kirchhoff potential $\phi'=a$ and establishes dissipation of $\int_\Omega\psi(u)\,\d\mathbf{x}$, where $\psi'=\phi$. The present work addresses general coefficients $a(t,\mathbf{x},u,\Grad u)$ and establishes unconditional $L^2$ stability without requiring a Kirchhoff potential. Even for $a=a(u)$, the two papers control different energies. Related resolvent expansions occur in kernel methods \cite{10.1007/s10915-020-01152-w,10.4208/cicp.OA-2019-0197}; in particular, the latter already contains a symmetric two-sided filter structure in the one-dimensional periodic setting with exact convolution. The comparison in \cref{sec:kernel} distinguishes this structure from the filtered coefficient evaluation and discrete stability theory developed here.

\subsection{Organization of the paper}

\Cref{sec:space} constructs diffusion discretizations for given coefficients, describes the compatible spatial realizations, and establishes the comparison principle used for the physical and splitting operators. \Cref{sec:criterion} motivates symmetric resolvents, develops high-order filters and SSP-EIN stability under stagewise splitting conditions, and then treats classical RK-based EIN methods and implementation. \Cref{sec:generalizations} treats continuous space, systems, and quadrature, and \cref{sec:kernel} compares the construction with kernel-based resolvent expansions.

\section{Discrete gradient discretization of diffusion operators}\label{sec:space}

Let $\mathcal T_h$ be a mesh partition of $\Omega$, and let $V_h \subset L^2(\Omega)$ and $\Sigma_h \subset [L^2(\Omega)]^d$ be finite-dimensional approximation spaces for the unknown and its discrete gradient. We take $\Sigma_h$ to be a broken finite element space with no interelement continuity constraints, allowing the gradient to be computed \textit{locally}. Functions in both spaces are assumed sufficiently regular on each element for the volume and trace pairings below. Write $(\cdot,\cdot)$ and $\norm{\cdot}$ for the exact scalar or vector $L^2$ product and norm, and $\Pi_{V_h}$ for the orthogonal projection onto $V_h$. The four Galerkin realizations in \cref{sec:spatial-realizations} use exact integration until \cref{sec:quadrature}; the SBP--SAT realization in \cref{sec:sbp-sat} uses prescribed discrete mass products.

\subsection{Discrete gradients and diffusion operators for given coefficients}

With $\mathbf n_K$ the outward unit normal of an element $K$, define $G_h: V_h\mapsto\Sigma_h$ by
\begin{equation}\label{eq:ldg-gradient}
	(G_h v_h, \mathbf r_h)_K = -(v_h, \Div\mathbf r_h)_K + \langle\widehat v_h, \mathbf r_h\cdot\mathbf n_K\rangle_{\partial K},
	\quad \forall \mathbf r_h \in\Sigma_h|_K.
\end{equation}
Here $\widehat v_h$ is a prescribed numerical trace, linear in $v_h$ and dependent only on immediately neighboring elements. Fix the spaces, gradient reconstruction, and stabilization structure. For a given nonnegative diffusion coefficient $\kappa(\mathbf{x})$ and nonnegative stabilization weights $\eta$, define
\begin{equation}\label{eq:generic-diffusion-form}
	\mathfrak a_h[\kappa; \eta](z_h, w_h) 
	\coloneqq \int_\Omega \kappa(\mathbf{x}) G_h z_h \cdot G_h w_h \, \d\mathbf{x} + c_h[\eta](z_h, w_h),
\end{equation}
where $c_h[\eta]$ is a symmetric semi-positive bilinear form. The corresponding discrete diffusion operator $\mathscr L_h[\kappa;\eta]: V_h\mapsto V_h$ is defined by
\begin{equation}\label{eq:generic-diffusion-operator}
	(\mathscr L_h[\kappa;\eta] z_h, w_h)
	\coloneqq \mathfrak a_h[\kappa;\eta](z_h, w_h),
	\quad\forall z_h,w_h\in V_h.
\end{equation}
It is well defined, self-adjoint, and semi-positive in the fixed mass product. We use $\mathscr L_h$ for this semi-positive diffusion operator and $L_h$ below for the signed evolution field.

\subsection{Galerkin realizations and interface stabilization}\label{sec:spatial-realizations}

The gradient form \eqref{eq:generic-diffusion-form} provides the common structure for the four Galerkin realizations \cite{10.1137/S0036142901384162,10.1016/j.cam.2015.10.024}.

For CG, take a \(H_0^1\)-conforming space \(V_h\), a broken gradient space \(\Sigma_h\) containing $\Grad V_h$, and $\widehat v_h$ equal to the trace of $v_h$; elementwise integration by parts then gives $G_h v_h = \Grad v_h$. For BR1, Brezzi, and LDG, use standard DG spaces \(V_h\) and \(\Sigma_h\). On an interior face $e = \partial K^+ \cap \partial K^-$ for two neighboring elements \(K^\pm\), let $\mathbf n^\pm$ be the corresponding outward normals and use the DG average and jump operators \cite{10.1137/S0036142901384162} 
\begin{equation*}
	\avg{v_h} \coloneqq \frac{1}{2}(v_h^+ + v_h^-),
	\quad
	\jump{v_h} \coloneqq v_h^+ \mathbf n^+ + v_h^- \mathbf n^-.
\end{equation*}
The general LDG numerical trace \cite{10.1137/S0036142997316712} is
\begin{equation}\label{eq:ldg-biased-trace}
	\widehat v_h = \avg{v_h} + \mathbf C_{12,e}\jump{v_h},
\end{equation}
where the facewise bias vector $\mathbf C_{12,e} = \mathcal{O}(1)$ is fixed independently of the diffusion and splitting coefficients. Taking $\mathbf C_{12,e} = 0$ gives the central lifted gradient used by BR1 and Brezzi. Alternating traces are special LDG choices: $\mathbf C_{12,e}\cdot\mathbf n^+ = \frac{1}{2}$ gives $\widehat v_h = v_h^+$, while $\mathbf C_{12,e}\cdot\mathbf n^+ = -\frac{1}{2}$ gives $\widehat v_h = v_h^-$. Such alternating choices are used in the minimal-dissipation LDG (MD-LDG) method \cite{10.1007/s10915-007-9130-3} on conforming simplicial meshes, together with zero interior flux penalty and its prescribed boundary treatment. Boundary traces incorporate the prescribed homogeneous boundary conditions, with periodic faces paired when periodic conditions are used.

The choices $c_h \equiv 0$ and the following two types cover the basic spatial stabilization mechanisms:
\begin{equation}\label{eq:penalties}
	\begin{aligned}
		c_h^{\mathrm{I}}[\eta](z_h, w_h)
		& \coloneqq \sum_e \int_e \eta_e^{\mathrm{I}} \jump{z_h} \cdot \jump{w_h} \, \d s,\\
		c_h^{\mathrm{II}}[\eta](z_h, w_h)
		& \coloneqq \sum_e \int_{\omega_e} \eta_e^{\mathrm{II}} L_e(\jump{z_h}) \cdot L_e(\jump{w_h}) \, \d\mathbf{x}.
	\end{aligned}
\end{equation}
Here $\eta$ collects the weights $\eta_e^{\mathrm{I}} \geq 0 $ or $\eta_e^{\mathrm{II}} \geq 0$, $L_e$ is a fixed face lifting into a patch space, and $\omega_e$ is its support. When this construction is applied to nonlinear diffusion, these weights may be evaluated from time, state, discrete gradient, and traces. On a Dirichlet boundary face, set $\jump{v_h}=v_h\mathbf n$, with $\mathbf n$ the outward unit normal. Sums of these forms are allowed. Type I is the face-jump stabilization used in penalized LDG methods. Type II is the local lifting stabilization used in the method of Brezzi et al. The choices are summarized in \cref{tab:spatial-realizations}. In every realization, $\kappa$ multiplies the fixed reconstructed volume gradient in \eqref{eq:generic-diffusion-form}. This is equivalent to projecting the complete constitutive flux onto $\Sigma_h$ before pairing with the test gradient. These forms allow $\kappa \geq 0$ and require no inverse coefficient.

\begin{table}[htbp]
	\centering
	\caption{CG, BR1, Brezzi, and LDG realizations and their spatial stabilization. Types I and II refer to \eqref{eq:penalties}.}
	\label{tab:spatial-realizations}
	\setlength{\tabcolsep}{3pt}
	\begin{tabular}{lll}
		\toprule
		Spatial & Discrete gradient \(G_h\) & Stabilization \(c_h\) \\
		\midrule
		CG & Strong gradient & Zero \\
		BR1 & Central lifted gradient & Zero \\
		Brezzi & Central lifted gradient & Type II \\
		Penalized LDG & $\mathbf C_{12}$-biased lifted gradient & Type I \\
		MD-LDG & Alternating lifted gradient & Type I with zero interior penalty \\
		\bottomrule
	\end{tabular}
\end{table}

\subsection{SBP--SAT realization}\label{sec:sbp-sat}

Summation-by-parts (SBP) discretizations with simultaneous approximation terms (SATs) fit the same gradient-form framework in a discrete mass norm; see \cite{10.1016/j.jcp.2021.110634} for diffusion SATs and their energy-stability and adjoint-consistency conditions. Let $H$ and $H_\Sigma$ be fixed symmetric positive definite state and gradient norm matrices, and let $G_h$ be a fixed discrete gradient incorporating any fixed boundary or interface lifting. For example, a one-dimensional SBP derivative $\mathsf D$ satisfies $H\mathsf D+\mathsf D^TH=B$, where $B$ is the signed boundary evaluation matrix.

For a given nonnegative coefficient $\kappa$, let $b_h[\kappa](\mathbf q,\mathbf r)$ be its symmetric semi-positive weighted form on the gradient space. For diagonal-norm nodal realizations,
\begin{equation}\label{eq:sbp-coefficient-form}
	b_h[\kappa](\mathbf q, \mathbf r) = \sum_\ell\omega_\ell \kappa(\mathbf x_\ell)
	\mathbf q_\ell\cdot\mathbf r_\ell,
\end{equation}
where $\omega_\ell>0$ are the gradient norm weights, repeated over spatial components in $H_\Sigma$. More general gradient metrics are allowed provided the coefficient form remains symmetric semi-positive. For homogeneous boundary data, require that the complete spatial discretization, including all SATs, admit the form
\begin{equation}\label{eq:sbp-factorization}
	\mathfrak a_h[\kappa;\eta](z,w)
	\coloneqq b_h[\kappa](G_hz,G_hw)+(W[\eta]Ez,Ew).
\end{equation}
Here $E$ is a fixed collection of boundary traces or interface residuals and $W[\eta]$ is symmetric semi-positive. The last term represents the remaining stabilization after consistency contributions have been incorporated into the gradient form. Define $\mathscr L_h[\kappa;\eta]$ by \eqref{eq:generic-diffusion-operator} using $(\cdot,\cdot)_H$; it is self-adjoint semi-positive in that product.

\subsection{Coefficient comparison and the physical and splitting operators}\label{sec:spatial-comparison}

\begin{lemma}[Comparison for given coefficients]\label{lem:spatial-comparison}
	For the Galerkin realizations \eqref{eq:generic-diffusion-form}, let $\kappa,s \geq 0$ and let $\eta, \eta_S$ be nonnegative stabilization weights for the same gradient and jump or lifting forms. If $\rho > 0$ satisfies $\rho\kappa\leq s$ and $\rho\eta\leq\eta_S$ pointwise on the respective volumes, faces, or patches, then
	\begin{equation}\label{eq:generic-majorant}
		0 \preceq\rho\mathscr L_h[\kappa;\eta]
		\preceq\mathscr L_h[s;\eta_S].
	\end{equation}
	For SBP--SAT realizations \eqref{eq:sbp-factorization}, the same conclusion holds in the state mass product if
	\begin{equation}\label{eq:sbp-form-comparison}
		\rho b_h[\kappa](\mathbf q,\mathbf q)\leq b_h[s](\mathbf q,\mathbf q)
		\quad\forall\mathbf q,\quad
		\rho W[\eta]\preceq W[\eta_S].
	\end{equation}
	For \eqref{eq:sbp-coefficient-form} and diagonal penalty weights, these conditions follow from the corresponding nodal coefficient and scalar weight bounds.
\end{lemma}

\begin{proof}
	Evaluate the forms at $(z_h, z_h)$ for any \(z_h \in V_h\). Bound each coefficient-weighted squared gradient, jump, or lifting term by its counterpart with coefficient $s$ and weights $\eta_S$, and sum up. The same argument applied to \eqref{eq:sbp-factorization} proves the SBP--SAT statement.
\end{proof}

We now apply this construction to both diffusion operators in \eqref{eq:split-intro}. For a fixed time $t$ and state $v_h$, set the frozen physical diffusion coefficient $a_{t,v_h}(\mathbf{x})\coloneqq a(t,\mathbf{x},v_h,G_hv_h)$. Let $\eta_A(t,v_h)$ collect the physical stabilization weights evaluated at this state, gradient, and traces, and let $\eta_S(t)$ be splitting weights. Define
\begin{equation}\label{eq:physical}
	\mathcal A_h(t,v_h) \coloneqq \mathscr L_h[a_{t,v_h};\eta_A(t,v_h)],
	\quad S_h(t) \coloneqq \mathscr L_h[s(t,\cdot);\eta_S(t)].
\end{equation}
For later use, denote the splitting bilinear form by
\begin{equation}\label{eq:splitting-bilinear-form}
	\mathfrak s_h(t)\coloneqq\mathfrak a_h[s(t,\cdot);\eta_S(t)].
\end{equation}
Both operators use the spatial construction above; the state dependence of $\mathcal A_h(t,v_h)$ produces the nonlinear evolution field. For SBP--SAT, $a_{t,v_h}$ is evaluated at the nodes, and $c_h$ denotes the remaining stabilization term in \eqref{eq:sbp-factorization}.

The Galerkin semidiscretization for \eqref{eq:model} is
\begin{equation}\label{eq:original-semidiscrete-scheme}
	(\partial_t u_h, w_h) + (\mathcal A_h(t, u_h) u_h, w_h) = 0,
	\quad\forall w_h \in V_h,
\end{equation}
with $u_h(0)=\Pi_{V_h} u_0$. Equivalently, one may write
\begin{equation}
	\partial_t u_h = L_h(t,u_h)\coloneqq-\mathcal A_h(t,u_h)u_h. 
\end{equation}
Testing with $u_h$ gives
\begin{equation}\label{eq:original-semidiscrete-stability}
	\frac{1}{2}\frac{\d}{\d t}\norm{u_h}^2 = -(\mathcal A_h(t, u_h) u_h, u_h) \leq 0,
\end{equation}
a discrete counterpart of \eqref{eq:semi-discrete-stability}. For SBP--SAT, the corresponding scheme and energy identity are
\begin{equation}\label{eq:sbp-scheme}
	(\partial_t u_h,w)_H + (\mathcal A_h(t, u_h) u_h, w)_H = 0, \quad\forall w,
\end{equation}
\begin{equation}\label{eq:sbp-energy}
	\frac{1}{2}\frac{\d}{\d t}\norm{u_h}_H^2 = -(\mathcal A_h(t,u_h) u_h, u_h)_H \leq 0.
\end{equation}
The later filtering and RK arguments apply in this fixed mass norm with the same sufficient splitting factors, as formalized in \cref{sec:quadrature}.

\section{Symmetric filtering and Runge--Kutta stability}\label{sec:criterion}

\subsection{Why a symmetric resolvent is useful}\label{sec:motivation}

Consider one step from $t^n$ to $t^{n+1} = t^n + \tau^n$, with input $u_h^n$ and a splitting operator $S_h^n = \mathscr L_h[s^n;\eta_S^n]$. Classical EIN1 treats the linear splitting term implicitly and the remainder explicitly \cite{10.1007/s11425-018-9524-x}:
\begin{equation*}
	\frac{u_h^{n+1} - u_h^n}{\tau^n}
	= -S_h^{n} u_h^{n+1} + S_h^{n} u_h^n - \mathcal A_h(t^n, u_h^n)u_h^n.
\end{equation*}
Equivalently,
\begin{equation}\label{eq:classical-ein1}
	u_h^{n+1} = u_h^{n} - \tau^n(I + \tau^n S_h^{n})^{-1} \mathcal{A}_h(t^n, u_h^{n}) u_h^{n}.
\end{equation}
Under $0 \preceq \frac{1}{2}\mathcal{A}_h(t^n,u_h^{n}) \preceq S_h^{n}$, \eqref{eq:classical-ein1} contracts the time-dependent norm induced by $P_h\coloneqq I + \tau^n S_h^{n}$ as in \cite[Theorem 3.1]{10.1007/s11425-018-9524-x}:
\begin{equation}\label{eq:modified}
	\norm{u_h^{n+1}}^2 + \tau^n (S_h^{n} u_h^{n+1}, u_h^{n+1}) 
	\leq \norm{u_h^{n}}^2 + \tau^n (S_h^{n} u_h^{n}, u_h^{n}),
\end{equation}
but this contraction cannot imply a discrete version of \eqref{eq:semi-discrete-stability}. Indeed, 
\begin{equation}
	P_h^{\frac{1}{2}}u_h^{n+1}
	= (I-\tau^n P_h^{-\frac{1}{2}} \mathcal{A}_h(t^n,u_h^{n})P_h^{-\frac{1}{2}}) P_h^{\frac{1}{2}} u_h^{n},
\end{equation}
where $0 \preceq \tau^n P_h^{-\frac{1}{2}} \mathcal A_h(t^n,u_h^{n}) P_h^{-\frac{1}{2}} \preceq 2I$. Hence the transformed update is a contraction. If $S_h^{n}$ and $\mathcal{A}_h(t^n, u_h^{n})$ commute, simultaneous diagonalization also proves contraction in the original norm. Without commutation the operator majorant alone is insufficient. For example,
\begin{equation}\label{eq:counterexample}
	S_h^{n} = \begin{pmatrix}1&0\\0&100\end{pmatrix},\quad
	\mathcal{A}_h = \begin{pmatrix}\frac{1}{2}&5\\5&50\end{pmatrix},\quad
	0 \preceq \mathcal{A}_h\preceq S_h^{n},
\end{equation}
but with $\tau^n=1$ and $u_h^{n}=(0,1)^T$, \eqref{eq:classical-ein1} gives $u_h^{n+1}=(-\frac{5}{2},\frac{51}{101})^T$ and $\norm{u_h^{n+1}} > \norm{u_h^{n}}$. This motivates a two-sided resolvent product.

\begin{proposition}[Two-sided resolvent stability]\label{prop:two-sided-ein1}
	Fix a self-adjoint semi-positive splitting operator $S_h^{n}$ and compute
	\begin{equation}\label{eq:two-sided-ein1}
		\begin{aligned}
			w_h & = (I + \frac{\tau^n}{2} S_h^{n})^{-1} u_h^n, \\
			u_h^{n+1} & = u_h^n - \tau^n (I + \frac{\tau^n}{2} S_h^{n})^{-1} \mathcal A_h(t^n, w_h) w_h.
		\end{aligned}
	\end{equation}
	Suppose, at this intermediate state \(w_h\), that for $\rho = \frac{1}{4}$,
	\begin{equation*}
		0 \leq \rho a(t^n, \mathbf x,w_h,G_hw_h) \leq s^n(\mathbf x),
		\quad 0\leq\rho\eta_A(t^n,w_h)\leq\eta_S^n.
	\end{equation*}
	The second inequality applies to each penalty weight. Then $\norm{u_h^{n+1}} \leq \norm{u_h^n}$ for every $\tau^n > 0$.
\end{proposition}

\begin{proof}
	By \cref{lem:spatial-comparison}, $0\preceq\rho\mathcal A_h(t^n,w_h)\preceq S_h^{n}$. The self-adjoint semi-positive operator
	\begin{equation*}
		Q_h = \tau^n (I + \frac{\tau^n}{2} S_h^{n})^{-1}\mathcal A_h(t^n,w_h)(I+\frac{\tau^n}{2}S_h^{n})^{-1}
	\end{equation*}
	satisfies $0\preceq Q_h\preceq\frac{1}{2\rho}I\preceq2I$, since $\frac{\tau^n\lambda}{(1+\frac{\tau^n}{2}\lambda)^2}\leq\frac{1}{2}$ for $\lambda\geq0$. Thus \eqref{eq:two-sided-ein1} is $u_h^{n+1} = (I-Q_h)u_h^n$, which contracts the fixed mass norm.
\end{proof}

\begin{remark}[An IMEX-type EIN update]\label{rem:two-sided-imex}
	Set $L_{I,h}^n(v_h)=-S_h^nv_h$ and $L_{E,h}^n(t,v_h)=S_h^nv_h-\mathcal A_h(t,v_h)v_h$. Then \eqref{eq:two-sided-ein1} becomes 
	\begin{equation*}
		\begin{aligned}
			w_h & = u_h^n + \frac{\tau^n}{2} L_{I,h}^{n}(w_h), \\
			u_h^{n+1} & = u_h^n + \tau^n L_{E,h}^{n}(t^n,w_h) + \tau^n L_{I,h}^{n}(w_h)
			+ \frac{\tau^n}{2}\left[L_{I,h}^{n}(u_h^{n+1})-L_{I,h}^{n}(u_h^n)\right].
		\end{aligned}
	\end{equation*}
	This gives an additive IMEX realization with explicit and implicit Butcher tableaux
	\begin{equation*}
		\begin{array}{c|ccc}
			0 & 0 & 0 & 0 \\
			0 & 0 & 0 & 0 \\
			1 & 0 & 1 & 0 \\
			\hline
			& 0 & 1 & 0
		\end{array},
		\quad
		\begin{array}{c|ccc}
			0 & 0 & 0 & 0 \\
			\frac{1}{2} & 0 & \frac{1}{2} & 0 \\
			1 & -\frac{1}{2} & 1 & \frac{1}{2} \\
			\hline
			& -\frac{1}{2} & 1 & \frac{1}{2}
		\end{array}.
	\end{equation*}
\end{remark}

For a time scale $\delta \geq 0$ and an integer $p \geq 1$, define the basic filter and the order-$p$ filter by 
\begin{equation}\label{eq:step-indexed-filters}
	J_{\delta,h}^n\coloneqq(I+\frac{\delta}{2}S_h^n)^{-1}, \quad
	R_{\delta,h}^n\coloneqq I-(I-J_{\delta,h}^n)^p.
\end{equation}
The latter one has the polynomial representation
\begin{equation}\label{eq:resolvent}
	R_{\delta,h}^n=\sum_{k=1}^p\alpha_k(J_{\delta,h}^n)^k,
	\quad\alpha_k\coloneqq(-1)^{k+1}\binom pk.
\end{equation}
Its scalar symbol and stability constant are defined by 
\begin{equation}\label{eq:sharp-ein-constant}
	r_p(x)\coloneqq1-\left(\frac{x}{1+x}\right)^p, \quad
	C_p^{\mathrm{EIN}} \coloneqq \max_{x \geq 0} xr_p(x)^2.
\end{equation}

\begin{lemma}[Sharp spectral bounds]\label{lem:sharp-resolvent}
	The maximum of \(x r_p(x)^2\) (\(x \geq 0\)) occurs at $x_p = \frac{y_p}{1-y_p}$, where $y_p\in(0,1)$ is the unique solution of
	\begin{equation}\label{eq:sharp-ein-root}
		2p y_p^p=\sum_{j=0}^{p-1}y_p^j,\quad
		C_p^{\mathrm{EIN}} = \frac{y_p}{1-y_p}(1-y_p^p)^2.
	\end{equation}
	In particular,
	\begin{equation}\label{eq:sharp-ein-values}
		\begin{gathered}
		C_1^{\mathrm{EIN}} = \frac{1}{4},\quad
		C_2^{\mathrm{EIN}} = \frac{107+51\sqrt{17}}{512}\approx0.6196843494,\\
		C_3^{\mathrm{EIN}}\approx1.0101037996,\quad
		C_4^{\mathrm{EIN}}\approx1.4078108284.
		\end{gathered}
	\end{equation}
	For every self-adjoint semi-positive operator $S_h^n$ and every $\delta>0$,
	\begin{equation}\label{eq:sharp-resolvent-operator}
		0\preceq\delta R_{\delta,h}^n S_h^n R_{\delta,h}^n \preceq 2 C_p^{\mathrm{EIN}} I.
	\end{equation}
	For $\rho>0$, whenever $0\preceq\rho\mathcal A_h(t,R_{\delta,h}^nv_h)\preceq S_h^n$, the operator ordering is
	\begin{equation}\label{eq:filtered-operator-chain}
		0\preceq\delta R_{\delta,h}^n\mathcal A_h(t,R_{\delta,h}^nv_h)R_{\delta,h}^n
		\preceq\frac{\delta}{\rho}R_{\delta,h}^nS_h^nR_{\delta,h}^n
		\preceq\frac{2C_p^{\mathrm{EIN}}}{\rho}I.
	\end{equation}
	Also, $C_p^{\mathrm{EIN}}<\frac{p^2}{4}$ for $p>1$.
\end{lemma}

\begin{proof}
	Set $y\coloneqq\frac{x}{1+x}$. Then $xr_p(x)^2=\frac{y}{1-y}(1-y^p)^2$, which extends continuously to zero at both endpoints of $[0,1]$. At an interior stationary point, differentiation gives $1-y^p=2p y^p(1-y)$, equivalently \eqref{eq:sharp-ein-root}. After division by $y^p$, the equation is $\sum_{j=1}^p y^{-j}=2p$. Its left side decreases strictly from infinity to $p$, proving uniqueness of the maximizer. For $p=2$, $y_2=\frac{1+\sqrt{17}}{8}$ gives the stated closed form. Functional calculus with $x=\frac{\delta}{2}S_h^n$ proves \eqref{eq:sharp-resolvent-operator}. Congruence by $R_{\delta,h}^n$ gives \eqref{eq:filtered-operator-chain}. Finally, the geometric sum gives $r_p(x)\leq\frac{p}{1+x}$, with strict inequality for $p>1$ and $x>0$. Evaluating at the finite positive maximizer yields $C_p^{\mathrm{EIN}}<\frac{p^2x_p}{(1+x_p)^2}\leq\frac{p^2}{4}$.
\end{proof}

At time $t$ and input $v_h$, define the filtered state, filtered Euler map and its dissipations by
\begin{equation}\label{eq:branch}
	\begin{aligned}
		w_h&\coloneqq R_{\delta,h}^nv_h,\\
		f_{\delta,h}^n(t,v_h)&\coloneqq-R_{\delta,h}^n\mathcal A_h(t,w_h)w_h,\\
		\mathcal T_{\delta,h}^n(t,v_h)&\coloneqq v_h+\delta f_{\delta,h}^n(t,v_h),
	\end{aligned}
\end{equation}
\begin{equation}\label{eq:dissipation}
	\mathcal D_{\delta,h}^n(t,v_h)\coloneqq(w_h,\mathcal A_h(t,w_h)w_h),\quad
	\mathcal L_{\delta,h}^n(t,v_h)\coloneqq\mathcal D_{\delta,h}^n(t,v_h)-\frac{\delta}{2}\norm{f_{\delta,h}^n(t,v_h)}^2.
\end{equation}
In the spatial realization of \cref{sec:space},
\begin{equation}\label{eq:physical-dissipation}
	\mathcal D_{\delta,h}^n(t,v_h)
	= \int_\Omega a(t,\mathbf x,w_h,G_hw_h)|G_hw_h|^2\,\d\mathbf x
	+ c_h[\eta_A(t,w_h)](w_h,w_h).
\end{equation}

\begin{theorem}[Stability and consistency of filtered forward Euler]\label{lem:forward}
	Let $S_h^n=\mathscr L_h[s^n;\eta_S^n]$ be self-adjoint and semi-positive. At $w_h=R_{\delta,h}^n v_h$, assume
	\begin{equation}\label{eq:euler-coefficient-bounds}
		0\leq\rho a(t,\mathbf x,w_h,G_hw_h)\leq s^n(\mathbf x),\quad
		0\leq\rho\eta_A(t,w_h)\leq\eta_S^n,
	\end{equation}
	for $\rho>0$, with the second inequality imposed on each penalty weight. Put $c_p=\frac{2 C_p^{\mathrm{EIN}}}{\rho}$ and $\eta_p=1-\frac{C_p^{\mathrm{EIN}}}{\rho}$. Then
	\begin{equation}\label{eq:forward-identity}
		\begin{aligned}
			\frac{1}{2}\norm{\mathcal T_{\delta,h}^n(t,v_h)}^2+\delta\mathcal L_{\delta,h}^n(t,v_h)&=\frac{1}{2}\norm{v_h}^2,\\
			\mathcal L_{\delta,h}^n(t,v_h)&\geq\eta_p\mathcal D_{\delta,h}^n(t,v_h).
		\end{aligned}
	\end{equation}
	Thus $\rho\geq C_p^{\mathrm{EIN}}$ implies $\norm{\mathcal T_{\delta,h}^n(t,v_h)}\leq\norm{v_h}$ for every $\delta>0$. Moreover,
	\begin{equation}\label{eq:branch-increment}
		\delta\mathcal L_{\delta,h}^n(t,v_h) \geq (c_p^{-1}-\frac{1}{2})
		\norm{\mathcal T_{\delta,h}^n(t,v_h) - v_h}^2.
	\end{equation}
	On a fixed spatial discretization, formally,
	\begin{equation}\label{eq:euler-consistency}
		\mathcal T_{\delta,h}^n(t,v_h) = v_h + \delta L_h(t,v_h) + \mathcal O(\delta^{p+1}).
	\end{equation}
	For $\delta=\tau^n$, the update formally differs from the original forward Euler step of \eqref{eq:original-semidiscrete-scheme} by $\mathcal O((\tau^n)^{p+1})$ and retains formal order one.
\end{theorem}

\begin{proof}
	By \cref{lem:spatial-comparison,lem:sharp-resolvent},
	\begin{equation}\label{eq:bounded-product}
		B_h\coloneqq\delta R_{\delta,h}^n\mathcal A_h(t,w_h)R_{\delta,h}^n,
		\quad 0\preceq B_h\preceq c_pI.
	\end{equation}
	Since $B_h^2 \preceq c_p B_h$,
	\begin{equation}\label{eq:stage-inequality}
		(v_h,f_{\delta,h}^n(t,v_h))=-\mathcal D_{\delta,h}^n(t,v_h),\quad
		\delta\norm{f_{\delta,h}^n(t,v_h)}^2\leq c_p\mathcal D_{\delta,h}^n(t,v_h).
	\end{equation}
	Expanding $\norm{v_h+\delta f_{\delta,h}^n(t,v_h)}^2$ proves the energy estimates. For formal consistency,
	\begin{equation}\label{eq:resolvent-consistency}
		I-J_{\delta,h}^n=\frac{\delta}{2}S_h^nJ_{\delta,h}^n=\mathcal O(\delta),\quad
		R_{\delta,h}^n-I=-(I-J_{\delta,h}^n)^p=\mathcal O(\delta^p).
	\end{equation}
	Formal expansion gives $f_{\delta,h}^n(t,v_h) - L_h(t,v_h)=\mathcal O(\delta^p)$ and hence \eqref{eq:euler-consistency}.
\end{proof}

\begin{remark}[Constant splitting coefficient]
	Taking $s^n(\mathbf x)\equiv s_0$ gives the constant-coefficient splitting used in \cite{10.1007/s11425-018-9524-x}. In the present filtered Euler method, a sufficient condition is
	\begin{equation*}
		s_0\geq C_p^{\mathrm{EIN}}a(t, \mathbf x, w_h, G_h w_h),
		\quad \eta_S^n\geq C_p^{\mathrm{EIN}}\eta_A(t,w_h),
	\end{equation*}
	at the filtered state $w_h = R_{\delta,h}^n v_h$, with the second bound imposed on each penalty weight. Then the filtered Euler update is $L^2$ contractive for every $\delta > 0$.
\end{remark}

\begin{remark}[Sharp Euler threshold]
	The condition $\rho\geq C_p^{\mathrm{EIN}}$ is sharp under the operator comparison: take $S_h^n = sI$, $\mathcal A_h=\frac{s}{\rho}I$, and $\frac{\delta s}{2}=x_p$. The amplification factor is $1-\frac{2C_p^{\mathrm{EIN}}}{\rho}$, whose absolute value exceeds one when $\rho < C_p^{\mathrm{EIN}}$.
\end{remark}

\subsection{SSP-EIN methods}\label{sec:filtered-euler}\label{sec:ssp-general}

Suppress $n$ on stage quantities in the following single-step formulas, retaining $\tau^n$ and the endpoint indices. Let an explicit SSP RK method of order $p$ have a Shu--Osher representation
\begin{equation}\label{eq:ssp-representation}
	\begin{aligned}
		& u_h^{(0)} = u_h^{n}, \\
		& u_h^{(i)} = \sum_{j<i}\alpha_{ij}\left[u_h^{(j)}+q_{ij}\tau^n L_h(t^{(j)},u_h^{(j)})\right],\quad 1\leq i\leq s,\\
		& u_h^{n+1} = u_h^{(s)},
	\end{aligned}
\end{equation}
where $\alpha_{ij}\geq0$, $\sum_{j<i}\alpha_{ij}=1$, and $q_{ij}\geq0$. Omit terms with $\alpha_{ij}=0$. The stage times are
\begin{equation}\label{eq:ssp-times}
	t^{(0)}\coloneqq t^n,\quad 
	t^{(i)}\coloneqq\sum_{j<i}\alpha_{ij}(t^{(j)}+q_{ij}\tau^n).
\end{equation}
At source stage $j$, choose a nonnegative splitting coefficient and penalty weights and set
\begin{equation}\label{eq:splitting-form}
	S_h^{n,(j)}\coloneqq\mathscr L_h[s^{n,(j)};\eta_S^{n,(j)}],\quad
	\mathfrak s_h^{n,(j)}\coloneqq\mathfrak a_h[s^{n,(j)};\eta_S^{n,(j)}].
\end{equation}
The corresponding EIN splitting operators are
\begin{equation}\label{eq:split}
	L_{E,h}^{n,(j)}(t,v_h)\coloneqq S_h^{n,(j)}v_h-\mathcal A_h(t,v_h)v_h,\quad
	L_{I,h}^{n,(j)}(v_h)\coloneqq-S_h^{n,(j)}v_h.
\end{equation}
Use \eqref{eq:step-indexed-filters}, \eqref{eq:branch}, and \eqref{eq:dissipation} with $S_h^{n,(j)}$ to define $J_{\delta,h}^{n,(j)}$, $R_{\delta,h}^{n,(j)}$, $f_{\delta,h}^{n,(j)}$, $\mathcal T_{\delta,h}^{n,(j)}$, and their dissipations.

For each branch, set 
\begin{equation}\label{eq:stage-filtered-branch-state}
	w_h^{n,(i,j)} \coloneqq R_{q_{ij}\tau^n,h}^{n,(j)}u_h^{n,(j)}.
\end{equation}
On branches with $\alpha_{ij} q_{ij}>0$, the coefficient and penalty comparisons, with $\rho^{n,(i,j)}>0$, are
\begin{equation}\label{eq:stage-substep-bounds}
	\begin{aligned}
		& 0\leq\rho^{n,(i,j)}a(t^{n,(j)},\mathbf x,w_h^{n,(i,j)},G_hw_h^{n,(i,j)})\leq s^{n,(j)}(\mathbf x),\\
		& 0\leq\rho^{n,(i,j)}\eta_e^\#(t^{n,(j)},w_h^{n,(i,j)})\leq\eta_{S,e}^{\#,n,(j)},\quad\#\in\{\mathrm I,\mathrm{II}\}.
	\end{aligned}
\end{equation}
Replace each Euler branch in \eqref{eq:ssp-representation} by its filtered map, then we obtain our SSP-EIN scheme. 
\begin{equation}\label{eq:ssp-ein}
	u_h^{(0)} = u_h^n, \quad 
	u_h^{(i)} = \sum_{j<i}\alpha_{ij}\mathcal T_{q_{ij}\tau^n,h}^{(j)}(t^{(j)},u_h^{(j)}) \quad 1\leq i \leq s \quad
	u_h^{n+1} = u_h^{(s)}.
\end{equation}

\begin{theorem}[Stability of SSP-EIN$p$]\label{thm:ssp-exact-general}
	Use order-$p$ filters in \eqref{eq:ssp-ein}. Suppose \eqref{eq:stage-substep-bounds} holds on every branch with $\alpha_{ij} q_{ij} > 0$. Define the branch data by
	\begin{equation}\label{eq:ssp-general-branches}
		\begin{aligned}
			w_h^{(i,j)}&\coloneqq R_{q_{ij}\tau^n,h}^{(j)}u_h^{(j)},\\
			f_h^{(i,j)}&\coloneqq -R_{q_{ij}\tau^n,h}^{(j)}\mathcal A_h(t^{(j)},w_h^{(i,j)})w_h^{(i,j)},\\
			y_h^{(i,j)}&\coloneqq u_h^{(j)} + q_{ij}\tau^n f_h^{(i,j)},\\
			\mathcal D_h^{(i,j)}&\coloneqq(w_h^{(i,j)},\mathcal A_h(t^{(j)},w_h^{(i,j)})w_h^{(i,j)}),\\
			\mathcal L_h^{(i,j)}&\coloneqq\mathcal D_h^{(i,j)}-\frac{q_{ij}\tau^n}2\norm{f_h^{(i,j)}}^2.
		\end{aligned}
	\end{equation}
	Thus $u_h^{(i)} = \sum_{j<i}\alpha_{ij} y_h^{(i,j)}$. Define the mixing dissipation 
	\begin{equation}\label{eq:ssp-general-mixing}
		\mathcal V_h^{(i)}\coloneqq\frac{1}{2}\sum_{0\leq j<k<i}\alpha_{ij}\alpha_{ik}\norm{y_h^{(i,j)}-y_h^{(i,k)}}^2,
	\end{equation}
	and the backward weights
	\begin{equation}\label{eq:ssp-backward-weights}
		\beta_s \coloneqq 1, \quad\beta_j\coloneqq\sum_{i=j+1}^s\beta_i\alpha_{ij},\quad j=s-1,\ldots,0.
	\end{equation}
	Then $\beta_j\geq0$, $\beta_0=1$, and there always holds the identity
	\begin{equation}\label{eq:ssp-general-exact}
		\frac{1}{2}\norm{u_h^{n+1}}^2
		+\tau^n\sum_{i=1}^s\beta_i\sum_{j<i}\alpha_{ij}q_{ij}\mathcal L_h^{(i,j)}
		+\sum_{i=1}^s\beta_i\mathcal V_h^{(i)}
		=\frac{1}{2}\norm{u_h^{n}}^2.
	\end{equation}
	Now suppose \(\rho^{(i,j)}\geq C_p^{\mathrm{EIN}}\). For positive-length branches set \(\eta_p^{(i,j)}\coloneqq 1-\frac{C_p^{\mathrm{EIN}}}{\rho^{(i,j)}}\geq0\). For zero-length branches set $\eta_p^{(i,j)}=0$. Then
	\begin{equation}\label{eq:ssp-general-dissipation}
		\frac{1}{2}\norm{u_h^{n+1}}^2
		+\tau^n\sum_{i=1}^s\beta_i\sum_{j<i}\alpha_{ij}q_{ij}\eta_p^{(i,j)}\mathcal D_h^{(i,j)}
		+\sum_{i=1}^s\beta_i\mathcal V_h^{(i)}
		\leq\frac{1}{2}\norm{u_h^{n}}^2.
	\end{equation}
	The stage norms satisfy
	\begin{equation}\label{eq:ssp-energy}
		\norm{u_h^{(i)}}\leq\norm{u_h^{n}}\quad(0\leq i\leq s),
		\quad \norm{u_h^{n+1}}\leq\norm{u_h^{n}},\quad\tau^n>0.
	\end{equation}
	The method has formal order $p$ on a fixed spatial discretization.
\end{theorem}

\begin{proof}
	Apply \cref{lem:forward} with scale $q_{ij}\tau^n$ and factor $\rho^{(i,j)}$ on each positive-length branch. Zero-length branches are identities. Thus $\mathcal L_h^{(i,j)}\geq\eta_p^{(i,j)}\mathcal D_h^{(i,j)}$ on the positive-length branches, and
	\begin{equation}\label{eq:ssp-branch-energy}
		\frac{1}{2}\norm{y_h^{(i,j)}}^2+q_{ij}\tau^n\mathcal L_h^{(i,j)}
		= \frac{1}{2}\norm{u_h^{(j)}}^2.
	\end{equation}
	The convex-combination identity gives
	\begin{equation}\label{eq:ssp-stage-energy}
		\frac{1}{2}\norm{u_h^{(i)}}^2+\tau^n\sum_{j<i}\alpha_{ij}q_{ij}\mathcal L_h^{(i,j)}+\mathcal V_h^{(i)}
		=\sum_{j<i}\frac{\alpha_{ij}}2\norm{u_h^{(j)}}^2.
	\end{equation}
	The backward recursion and the row sums of $\alpha$ give $\beta_j\geq0$ and $\beta_0=\beta_s=1$. Multiplication by $\beta_i$ and summation cancel the intermediate stage energies and yield \eqref{eq:ssp-general-exact}. Nonnegativity and convexity give \eqref{eq:ssp-general-dissipation} and \eqref{eq:ssp-energy}.

	By \eqref{eq:euler-consistency}, filtering formally perturbs each Euler branch by $\mathcal O((\tau^n)^{p+1})$. Induction through the explicit recursion, with the stage times \eqref{eq:ssp-times}, gives the same formal difference from the unfiltered RK step, preserving order $p$.
\end{proof}

\begin{remark}[Constant splitting coefficient]
	The SSP construction also permits the constant splitting choice of \cite{10.1007/s11425-018-9524-x}: set $s^{(j)}(\mathbf x)\equiv s_0$ at every stage. Define the physical coefficient at each filtered branch by
	\begin{equation*}
		a_h^{(i,j)}(\mathbf x)
		\coloneqq a(t^{(j)},\mathbf x,w_h^{(i,j)},G_hw_h^{(i,j)}).
	\end{equation*}
	If $s_0\geq C_p^{\mathrm{EIN}}a_h^{(i,j)}(\mathbf x)$ throughout the domain on every branch with $\alpha_{ij}q_{ij}>0$, and the penalty-weight bounds hold with the same factor, then $\norm{u_h^{n+1}}\leq\norm{u_h^n}$ for every $\tau^n>0$. Thus a single constant dominating all filtered-stage diffusion coefficients suffices for the full SSP step.
\end{remark}

\begin{example}[SSP-EIN1]\label{ex:ein1}
	Use $p=1$, $R_{\tau^n,h}^{(0)} = J_{\tau^n,h}^{(0)}$, and
	\begin{equation}\label{eq:ein1}
		u_h^{n+1} = \mathcal T_{\tau^n,h}^{(0)}(t^n,u_h^n).
	\end{equation}
	This is \eqref{eq:two-sided-ein1}. The penalty factor is $C_1^{\mathrm{EIN}} = \frac{1}{4}$.
\end{example}

\begin{example}[SSP-EIN2]\label{ex:ein2}
	Use $p=2$, $R_{\tau^n,h}^{(i)} = 2 J_{\tau^n,h}^{(i)}-(J_{\tau^n,h}^{(i)})^2$, and the Heun formula
	\begin{equation}\label{eq:ein2}
		\begin{aligned}
			u_h^{(1)}&=\mathcal T_{\tau^n,h}^{(0)}(t^n,u_h^n),\\
			u_h^{n+1}&=\frac{1}{2}u_h^n+\frac{1}{2}\mathcal T_{\tau^n,h}^{(1)}(t^{n+1},u_h^{(1)}).
		\end{aligned}
	\end{equation}
	At both evaluation stages require the penalty factor $C_2^{\mathrm{EIN}}\approx0.6196843494$.
\end{example}

\begin{example}[SSP-EIN3]\label{ex:ein3}
	Use $p=3$, $R_{\tau^n,h}^{(i)}=3J_{\tau^n,h}^{(i)}-3(J_{\tau^n,h}^{(i)})^2+(J_{\tau^n,h}^{(i)})^3$, and
	\begin{equation}\label{eq:ein3}
		\begin{aligned}
			u_h^{(1)}&=\mathcal T_{\tau^n,h}^{(0)}(t^n,u_h^n),\\
			u_h^{(2)}&=\frac34u_h^n+\frac{1}{4}\mathcal T_{\tau^n,h}^{(1)}(t^{n+1},u_h^{(1)}),\\
			u_h^{n+1}&=\frac13u_h^n+\frac23\mathcal T_{\tau^n,h}^{(2)}(t^n+\frac{\tau^n}{2},u_h^{(2)}).
		\end{aligned}
	\end{equation}
	At all three evaluation stages require the penalty factor $C_3^{\mathrm{EIN}}\approx1.0101037996$.
\end{example}

In the above examples, every positive-length Euler branch has $q_{ij}=1$, so its filter time scale is $\tau^n$. The underlying forward Euler, Heun, and SSP-RK3 formulas have SSP coefficient one \cite{10.1137/S003614450036757X}. Applying \cref{thm:ssp-exact-general} gives unconditional $L^2$ stability and formal temporal orders one, two, and three.

\begin{example}[SSP-EIN(5,4)]\label{sec:ein4}
	Use the five-stage, fourth-order SSPRK(5,4) coefficients of Spiteri and Ruuth in the Shu--Osher form reported in \cite{10.1007/s10915-004-4635-5}. With order-four filters,
	\begin{equation}\label{eq:ein4}
		u_h^{(i)}=\sum_{j<i}\alpha_{ij}\mathcal T_{q_{ij}\tau^n,h}^{(j)}(t^{(j)},u_h^{(j)}),\quad i=1,\ldots,5,
	\end{equation}
	where $u_h^{(0)}=u_h^n$ and $u_h^{n+1}=u_h^{(5)}$, and the times satisfy \eqref{eq:ssp-times}. The branchwise conditions \eqref{eq:stage-substep-bounds} with $\rho^{(i,j)}\geq C_4^{\mathrm{EIN}}$ give unconditional $L^2$ stability by \cref{thm:ssp-exact-general}. The method has formal order four.
\end{example}

\begin{remark}[Stagewise splitting and additive IMEX-RK interpretation]\label{rem:ssp-additive}
	The basic-filter solves are linear implicit solves with $L_{I,h}^{(j)}$, while the physical coefficient is evaluated explicitly at a filtered state. When the stage operators coincide, auxiliary stages for these solves give an additive IMEX-RK realization of the SSP construction. The explicit and implicit row sums may differ \cite{10.1007/s10915-004-4636-4}. Stage-dependent operators instead give a family of linear implicit fields indexed by the physical stage times $t^{(j)}$. For nonautonomous diffusion, these times are those of the outer explicit RK method, as in \eqref{eq:ssp-times}.
\end{remark}

\subsection{Classical RK based EIN methods}\label{sec:outer-criterion}

Let $(A,b)$ be an explicit RK method of order $p$, with $c_i=\sum_{j<i}a_{ij}$ and $t^{(i)}=t^n+c_i\tau^n$. The original classical RK time discretization is 
\begin{equation}\label{eq:unfiltered-rk}
	\begin{aligned}
		& v_h^{(i)}=u_h^n+\tau^n\sum_{j<i}a_{ij}L_h(t^{(j)},v_h^{(j)}),\\
		& v_h^{n+1}=u_h^n+\tau^n\sum_i b_iL_h(t^{(i)},v_h^{(i)}).
	\end{aligned}
\end{equation}
Choose stage splitting operators by \eqref{eq:splitting-form} and use the common full-step scale $\delta=\tau^n$. The filtered RK formula is
\begin{equation}\label{eq:outer-rk}
	\begin{aligned}
		u_h^{(i)}&=u_h^n+\tau^n\sum_{j<i}a_{ij}f_{\tau^n,h}^{(j)},\\
		u_h^{n+1}&=u_h^n+\tau^n\sum_i b_i f_{\tau^n,h}^{(i)},
	\end{aligned}
\end{equation}
where 
\begin{equation}\label{eq:stage-filtered-residual}
	w_h^{(i)}\coloneqq R_{\tau^n,h}^{(i)}u_h^{(i)},\quad
	f_{\tau^n,h}^{(i)}\coloneqq-R_{\tau^n,h}^{(i)}\mathcal A_h(t^{(i)},w_h^{(i)})w_h^{(i)}.
\end{equation}
Write
\begin{equation}\label{eq:coefficient-evaluation}
	a_h^{(i)}(\mathbf x)\coloneqq a(t^{(i)},\mathbf x,w_h^{(i)},G_hw_h^{(i)}).
\end{equation}

For $b_i>0$, define
\begin{equation}\label{eq:outer-matrices}
	B\coloneqq\diag(b_i),\quad M\coloneqq BA+A^TB-bb^T,\quad C_\gamma\coloneqq M+\gamma B
\end{equation}
and the RK factor
\begin{equation}\label{eq:positive-rk-factor}
	\gamma_{\mathrm{RK}}\coloneqq-\lambda_{\min}(B^{-\frac{1}{2}}MB^{-\frac{1}{2}})>0.
\end{equation}

\begin{theorem}[Stability and order of filtered RK methods]\label{thm:outer}\label{cor:positive-rk}
	Let RK method $(A,b)$ have order $p$ and positive weights, and use order-$p$ full-step filters in \eqref{eq:outer-rk}. At every filtered stage, suppose
	\begin{equation}\label{eq:pointwise-majorant}
		\begin{aligned}
			& 0\leq\rho^{(i)}a_h^{(i)}(\mathbf x)\leq s^{(i)}(\mathbf x),\\
			& 0\leq\rho^{(i)}\eta_e^\#(t^{(i)},w_h^{(i)})\leq\eta_{S,e}^{\#,(i)},\quad\#\in\{\mathrm I,\mathrm{II}\},
		\end{aligned}
	\end{equation}
	with
	\begin{equation}\label{eq:positive-rk-bound}
		\rho^{(i)}\geq\gamma_{\mathrm{RK}}C_p^{\mathrm{EIN}}.
	\end{equation}
	For SBP--SAT, impose the corresponding form comparison \eqref{eq:sbp-form-comparison}. Then $\norm{u_h^{n+1}}\leq\norm{u_h^n}$ for every $\tau^n>0$. More generally, for any $\gamma\geq0$ with $C_\gamma\succeq0$ and $\rho^{(i)} \geq \gamma C_p^{\mathrm{EIN}}$, set $c_p^{(i)}=\frac{2C_p^{\mathrm{EIN}}}{\rho^{(i)}}$. Then
	\begin{equation}\label{eq:outer-estimate}
		\begin{aligned}
			\frac{1}{2}\norm{u_h^{n+1}}^2
			&+\tau^n\sum_i b_i(1-\frac{\gamma c_p^{(i)}}{2})\mathcal D_{\tau^n,h}^{(i)}(t^{(i)},u_h^{(i)})\\
			&+\frac{(\tau^n)^{2}}{2}\sum_{i,j}(C_\gamma)_{ij}(f_{\tau^n,h}^{(i)},f_{\tau^n,h}^{(j)})
			\leq\frac{1}{2}\norm{u_h^{n}}^2.
		\end{aligned}
	\end{equation}
	The method has formal order $p$.
\end{theorem}

\begin{proof}
	By \cref{lem:spatial-comparison}, the coefficient and penalty bounds imply
	\begin{equation}\label{eq:majorant}
		0\preceq\rho^{(i)}\mathcal A_h(t^{(i)},w_h^{(i)})\preceq S_h^{(i)}.
	\end{equation}
	Since $M_{ii}=-b_i^2$, the matrix $B^{-\frac{1}{2}}M B^{-\frac{1}{2}}$ has a negative smallest eigenvalue. Thus $\gamma_{\mathrm{RK}}>0$ and $C_{\gamma_{\mathrm{RK}}}\succeq0$.
	The RK squared-norm identity \cite{10.1137/18M122892X,10.1137/22M1472218} gives
	\begin{equation}\label{eq:rk-identity}
		\begin{aligned}
			\frac{1}{2}\norm{u_h^{n+1}}^2
			&+\tau^n\sum_i b_i\left[\mathcal D_{\tau^n,h}^{(i)}(t^{(i)},u_h^{(i)})-\frac{\gamma\tau^n}{2}\norm{f_{\tau^n,h}^{(i)}}^2\right]\\
			&+\frac{(\tau^n)^2}{2}\sum_{i,j}(C_\gamma)_{ij}(f_{\tau^n,h}^{(i)},f_{\tau^n,h}^{(j)})
			=\frac{1}{2}\norm{u_h^n}^2.
		\end{aligned}
	\end{equation}
	Apply \eqref{eq:stage-inequality} with $\delta=\tau^n$ at each stage. The last sum is nonnegative because $C_\gamma\succeq0$.

	By \eqref{eq:resolvent-consistency}, formally $f_{\tau^n,h}^{(i)}(t,v_h)-L_h(t,v_h)=\mathcal O((\tau^n)^p)$. Induction through the explicit recursion gives an $\mathcal O((\tau^n)^{p+1})$ formal difference between the filtered and unfiltered RK steps, preserving order $p$.
\end{proof}

\begin{remark}[Constant splitting coefficient]
	With the constant splitting choice $s^{(i)}(\mathbf x)\equiv s_0$ of \cite{10.1007/s11425-018-9524-x}, the present RK construction is unconditionally $L^2$ stable provided \(s_0\geq\gamma_{\mathrm{RK}}C_p^{\mathrm{EIN}}a_h^{(i)}(\mathbf x)\) at every filtered stage and the penalty-weight bounds hold with the same factor. 
\end{remark}

\begin{example}[RK4-EIN]\label{sec:rk4-ein}
	At stages $i=0,1,2,3$, use the order-four filters $R_{\tau^n,h}^{(i)}$ with independently selected $S_h^{(i)}$ and the stage residuals $f_{\tau^n,h}^{(i)}$ from \eqref{eq:stage-filtered-residual}, at times $(t^{(0)},t^{(1)},t^{(2)},t^{(3)})=(t^n,t^n+\frac{\tau^n}{2},t^n+\frac{\tau^n}{2},t^{n+1})$. Then
	\begin{equation}\label{eq:rk4-ein}
		\begin{aligned}
			u_h^{(0)}&\coloneqq u_h^{n},\\
			u_h^{(1)}&\coloneqq u_h^{n}+\frac{\tau^n}{2} f_{\tau^n,h}^{(0)},\\
			u_h^{(2)}&\coloneqq u_h^{n}+\frac{\tau^n}{2} f_{\tau^n,h}^{(1)},\\
			u_h^{(3)}&\coloneqq u_h^{n}+\tau^n f_{\tau^n,h}^{(2)},\\
			u_h^{n+1}&=u_h^{n}+\frac{\tau^n}{6}(f_{\tau^n,h}^{(0)}+2f_{\tau^n,h}^{(1)}+2f_{\tau^n,h}^{(2)}+f_{\tau^n,h}^{(3)}).
		\end{aligned}
	\end{equation}
	Suppose the spatial assumptions hold at every filtered RK4 stage, and choose
	\begin{equation}\label{eq:rk4-ein-bound}
		0\leq C_4^{\mathrm{EIN}}a_h^{(i)}(\mathbf{x})\leq s^{(i)}(\mathbf{x}),\quad i=0,1,2,3,
	\end{equation}
	with the same factor for splitting penalty weights at each stage. Then \eqref{eq:rk4-ein} satisfies $\norm{u_h^{n+1}}\leq\norm{u_h^{n}}$ for every $\tau^n>0$ and has formal order four by \cref{cor:positive-rk}.

	For the classical RK4 arrays,
	the matrix with $\gamma=1$ is
	\begin{equation}\label{eq:rk4-certificate}
		C_1=BA+A^TB-bb^T+B
		=\frac1{36}\begin{pmatrix}
			5&4&-2&-1\\4&8&2&-2\\-2&2&8&4\\-1&-2&4&5
		\end{pmatrix}\succeq0.
	\end{equation}
	The eigenvalues of $C_1$ are $0$, $\frac{1}{3}$, and $\frac{7\pm\sqrt{13}}{36}$. Since $\rho^{(i)}\geq C_4^{\mathrm{EIN}}$ at each stage, \cref{thm:outer} applies. 
\end{example}

\section{Numerical test}

\subsection{Implementation}\label{sec:implementation}

\begin{table}[htbp]
	\centering
	\caption{Sufficient stagewise splitting factors for $L^2$ stability. SSP filters use the branch scale $q_{ij}\tau^n$; classical RK filters use $\tau^n$.}
	\label{tab:summary}
	\begin{tabular}{ll}
		\toprule
		Method or construction & Factor $\rho_{\min}$ \\
		\midrule
		SSP-EIN$p$ & $C_p^{\mathrm{EIN}}$ \\
		Positive-weight RK & $\gamma_{\mathrm{RK}} C_p^{\mathrm{EIN}}$\\
		\midrule
		SSP-EIN1 & $C_1^{\mathrm{EIN}} = \frac{1}{4}$\\
		SSP-EIN2 & $C_2^{\mathrm{EIN}} \approx 0.619684$\\
		SSP-EIN3 & $C_3^{\mathrm{EIN}} \approx 1.010104$\\
		SSP-EIN(5,4) & $C_4^{\mathrm{EIN}} \approx 1.407811$\\
		RK4-EIN & $C_4^{\mathrm{EIN}} \approx 1.407811$\\
		\bottomrule
	\end{tabular}
\end{table}

\begin{algorithm}[htbp]
	\caption{One step of SSP-EIN$p$}\label{alg:ssp-ein-step}
	\begin{algorithmic}[1]
		\Require $u_h^n$, $t^n$, $\tau^n>0$, and the Shu--Osher coefficients $\alpha_{ij}$ and $q_{ij}$ of an $s$-stage SSP method.
		\Ensure $u_h^{n+1}$.
		\State Set $u_h^{(0)}=u_h^n$ and compute the stage times $t^{(j)}$ from \eqref{eq:ssp-times}.
		\For{$j=0,\ldots,s-1$}
			\State Choose $s^{(j)}$ and $\eta_S^{(j)}$, and form $S_h^{(j)}=\mathscr L_h[s^{(j)};\eta_S^{(j)}]$.
			\For{each $i>j$ with $\alpha_{ij}>0$}
				\If{$q_{ij}>0$}
					\State $w_h^{(i,j)}=R_{q_{ij}\tau^n,h}^{(j)}u_h^{(j)}$.
					\State $y_h^{(i,j)}=u_h^{(j)}-q_{ij}\tau^n R_{q_{ij}\tau^n,h}^{(j)}
					\mathcal A_h(t^{(j)},w_h^{(i,j)})w_h^{(i,j)}$.
				\Else
					\State $y_h^{(i,j)}=u_h^{(j)}$.
				\EndIf
				\State Store $y_h^{(i,j)}$ for use in subsequent stages.
			\EndFor
			\State Form the next stage from the stored updates:
			\[
				u_h^{(j+1)}=\sum_{\substack{0\leq k\leq j\\\alpha_{j+1,k}>0}}
				\alpha_{j+1,k}y_h^{(j+1,k)}.
			\]
		\EndFor
		\State $u_h^{n+1}=u_h^{(s)}$.
	\end{algorithmic}
\end{algorithm}

\begin{algorithm}[htbp]
	\caption{One step of classical RK-EIN}\label{alg:rk-ein-step}
	\begin{algorithmic}[1]
		\Require $u_h^n$, $t^n$, $\tau^n>0$, and the Butcher coefficients $a_{ij}$, $b_i$, and $c_i$ of an explicit $s$-stage RK method.
		\Ensure $u_h^{n+1}$.
		\For{$i=1,\ldots,s$}
			\State Compute the stage value and time:
			\[
				u_h^{(i)}=u_h^n+\tau^n\sum_{j=1}^{i-1}a_{ij}f_{\tau^n,h}^{(j)},
				\qquad t^{(i)}=t^n+c_i\tau^n.
			\]
			\State Choose $s^{(i)}$ and $\eta_S^{(i)}$, and form $S_h^{(i)}=\mathscr L_h[s^{(i)};\eta_S^{(i)}]$.
			\State $w_h^{(i)}=R_{\tau^n,h}^{(i)}u_h^{(i)}$.
			\State $f_{\tau^n,h}^{(i)}=-R_{\tau^n,h}^{(i)}
			\mathcal A_h(t^{(i)},w_h^{(i)})w_h^{(i)}$.
			\State Store $f_{\tau^n,h}^{(i)}$ for use in subsequent stages and the final update.
		\EndFor
		\State $u_h^{n+1}=u_h^n+\tau^n\sum_{i=1}^s b_i f_{\tau^n,h}^{(i)}$.
	\end{algorithmic}
\end{algorithm}

\Cref{alg:ssp-ein-step,alg:rk-ein-step} advance the solution from $t^n$ to $t^{n+1}=t^n+\tau^n$ using the order-$p$ filters in \eqref{eq:step-indexed-filters} and SSP or classical RK compositions. Each filtered residual requires two filter applications and one evaluation of the nonlinear diffusion operator. Since each filter uses $p$ basic-resolvent solves, the cost is $2p$ linear solves per residual. Without stage restarts, SSP-EIN1, SSP-EIN2, SSP-EIN3, and RK4-EIN therefore require respectively 1, 2, 3, and 4 nonlinear evaluations and 2, 8, 18, and 32 linear solves per step. For SSP-EIN(5,4), each distinct filtered residual requires one nonlinear evaluation and eight linear solves.

The basic solves involve $I+\frac{\delta}{2}S_h^{(j)}$, which is positive definite even when $S_h^{(j)}$ is only semi-positive definite. In matrix form, these are symmetric positive definite systems with the discrete mass matrix included. A factorization or preconditioner can be reused whenever the splitting operator and filter scale coincide. In \cref{alg:ssp-ein-step}, branches from the same stage with the same $q_{ij}$ also share the filtered state and residual, which need only be computed once.

\begin{remark}[Checking stage bounds]\label{rem:stage-restart}
	One should choose trial splitting coefficients and penalty weights before each stage, using the factor in \cref{tab:summary} as lower bounds. After computing the filtered states, check \eqref{eq:stage-substep-bounds} for SSP-EIN or \eqref{eq:pointwise-majorant} for classical RK-EIN. If a bound fails, one may increase the splitting coefficient and weights and restart that stage. In \cref{alg:ssp-ein-step}, recompute all filtered states $w_h^{(i,j)}$ and Euler updates $y_h^{(i,j)}$ associated with the revised $S_h^{(j)}$ before forming $u_h^{(j+1)}$. In \cref{alg:rk-ein-step}, recompute $w_h^{(i)}$ and $f_{\tau^n,h}^{(i)}$ before advancing to the next stage. This strategy is similar to \cite[Section 2.5]{10.1007/s42967-023-00359-6}. If the diffusion coefficient \(a\) is itself bounded, then a simple choice of constant splitting coefficient \(s \geq C_p^\mathrm{EIN} \sup a\) for SSP-EIN or \(s \geq \gamma_\mathrm{RK} C_p^\mathrm{EIN} \sup a\) for classical RK-EIN suffices.
\end{remark}

\subsection{Manufactured solution test}

We verify the accuracy of our proposed two-sided filtered SSP-EIN3 in \cref{alg:ssp-ein-step} and RK4-EIN \cref{alg:rk-ein-step} using a one-dimensional manufactured solution \cite[Example 5.1]{10.1007/s11425-018-9524-x}. Consider
\begin{equation*}
	u_t = \partial_x(a(x,u) \partial_x u) + f(t,x),
	\quad a(x,u) = \sin^2(u),
\end{equation*}
on $[-\pi,\pi]$ with periodic boundary conditions. We prescribe the exact solution $u(t,x)=\sin(x-t)$. The diffusion coefficient is nonnegative and bounded by \(1\), and vanishes where $u=0$. The source term is evaluated explicitly at each RK stage and combined with the nonlinear diffusion residual before applying the outer filter. We use uniform grids with \(N\) elements, use \(h\) to denote the mesh size, and choose \(\Delta t = 0.2 h\) with the last step shortened if necessary. The final time is taken to be \(T = 0.5\).

For simplicity, we use \(a_{\max} = 1\) for estimate of the upper bound of diffusion coefficient and take the constant splitting coefficient \(s \equiv 1.2 a_{\max} C_p^\mathrm{EIN}\). The convergence history is listed in \cref{tab:1d-ein3,tab:1d-ein4}.

\begin{table}[htbp]
	\centering
	\caption{One-dimensional manufactured solution at $T=0.5$: SSP-EIN3 with CG/LDG-$P_2$, $\Delta t=0.2h$, $h=2\pi/N$. All volume integrals and $L^2$ errors use 3-point Gauss quadrature; the last step is shortened to reach $T$.}
	\label{tab:1d-ein3}
	\begin{tabular}{rcccc}
		\hline
		$N$ & CG $L^2$ error & Order & LDG $L^2$ error & Order \\
		\hline
		80 & $4.804464\times10^{-6}$ & -- & $1.770630\times10^{-5}$ & -- \\
		160 & $6.018329\times10^{-7}$ & 2.997 & $2.231552\times10^{-6}$ & 2.988 \\
		320 & $7.525203\times10^{-8}$ & 3.000 & $2.780758\times10^{-7}$ & 3.004 \\
		640 & $9.430260\times10^{-9}$ & 2.996 & $3.526204\times10^{-8}$ & 2.979 \\
		1280 & $1.178782\times10^{-9}$ & 3.000 & $4.401522\times10^{-9}$ & 3.002 \\
		\hline
	\end{tabular}
\end{table}

\begin{table}[htbp]
	\centering
	\caption{One-dimensional manufactured solution at $T=0.5$: RK4-EIN with CG/LDG-$P_3$, $\Delta t=0.2h$, $h=2\pi/N$. All volume integrals and $L^2$ errors use 4-point Gauss quadrature; the last step is shortened to reach $T$.}
	\label{tab:1d-ein4}
	\begin{tabular}{rcccc}
		\hline
		$N$ & CG $L^2$ error & Order & LDG $L^2$ error & Order \\
		\hline
		80 & $3.940176\times10^{-8}$ & -- & $2.290565\times10^{-7}$ & -- \\
		160 & $2.523014\times10^{-9}$ & 3.965 & $1.438353\times10^{-8}$ & 3.993 \\
		320 & $1.599392\times10^{-10}$ & 3.980 & $8.976301\times10^{-10}$ & 4.002 \\
		640 & $1.006326\times10^{-11}$ & 3.990 & $5.652984\times10^{-11}$ & 3.989 \\
		1280 & $6.353682\times10^{-13}$ & 3.985 & $3.544970\times10^{-12}$ & 3.995 \\
		\hline
	\end{tabular}
\end{table}

\section{Generalizations}\label{sec:generalizations}

\subsection{Continuous spatial operators}\label{sec:continuous}

The stability argument can also be formulated before discretizing space. Set $X\coloneqq L^2(\Omega)$ and initially take admissible states in $V\coloneqq H_0^1(\Omega)$. On an abstract step, define
\begin{equation}\label{eq:continuous-forms}
	\begin{aligned}
		\mathfrak a_{t,v}(z,w)&\coloneqq\int_\Omega a(t,\mathbf{x},v,\Grad v)\Grad z\cdot\Grad w\,\d\mathbf{x},\\
		\mathfrak s^{(i)}(z,w)&\coloneqq\int_\Omega s^{(i)}(\mathbf{x})\Grad z\cdot\Grad w\,\d\mathbf{x}.
	\end{aligned}
\end{equation}
For classical RK methods, assume at each stage $s^{(i)}\geq0$, $\rho^{(i)}>0$, and
$0\leq\rho^{(i)}a(t^{(i)},\mathbf{x},w^{(i)},\Grad w^{(i)})\leq s^{(i)}(\mathbf{x})$
at $w^{(i)}=R_{\tau^n}^{(i)}u^{(i)}$. For SSP methods, use the continuous branch states $w^{(i,j)}=R_{q_{ij}\tau^n}^{(j)}u^{(j)}$ and branchwise factors $\rho^{(i,j)}$ as in \eqref{eq:stage-substep-bounds}. The coefficient comparison gives
$0\leq\rho^{(i)}\mathfrak a_{t^{(i)},w^{(i)}}(z,z)\leq\mathfrak s^{(i)}(z,z)$.
Assume the splitting forms admit semi-positive self-adjoint realizations $S^{(i)}$ in $X$, and all required weak solves and filtered coefficient evaluations are well defined in a common admissible form space. For the classical RK full-step filter, the weak resolvent is characterized, whenever the right-hand side is admissible, by
\begin{equation}\label{eq:weak-resolvent}
	(J_{\tau^n}^{(i)}\ell,z)+\frac{\tau^n}{2}\mathfrak s^{(i)}(J_{\tau^n}^{(i)}\ell,z)=\langle\ell,z\rangle.
\end{equation}
Use the continuous analogues of \eqref{eq:resolvent} and \eqref{eq:branch}, with stagewise $J_{\tau^n}^{(i)}$ and $R_{\tau^n}^{(i)}$ and weak dual right-hand sides. For SSP methods replace $\tau^n$ by the corresponding branch length.

\begin{theorem}[Continuous spatial stability]\label{thm:continuous}
	Under these form and admissibility assumptions, the Euler identity \eqref{eq:forward-identity} and the general SSP identity \eqref{eq:ssp-general-exact} hold in $X$, with
	$\mathcal D_{\tau^n}^{(i)}(t,v)=\mathfrak a_{t,R_{\tau^n}^{(i)}v}(R_{\tau^n}^{(i)}v,R_{\tau^n}^{(i)}v)$.
	For EIN1--3, the sufficient condition is $\rho^{(i)}\geq C_p^{\mathrm{EIN}}$ at every evaluation stage. For general SSP methods, impose the branchwise bounds \eqref{eq:stage-substep-bounds} with $\rho^{(i,j)}\geq C_p^{\mathrm{EIN}}$; classical RK methods use the factor in \eqref{eq:positive-rk-bound}. After restoring step indices, the estimates sum with only $\frac{1}{2}\norm{u_0}^2$ on the right.
\end{theorem}

\begin{proof}
	At each fixed stage the spectral estimate gives
	\begin{equation*}
		\tau^n\mathfrak s^{(i)}(R_{\tau^n}^{(i)}y,R_{\tau^n}^{(i)}y)\leq2C_p^{\mathrm{EIN}}\norm{y}^2.
	\end{equation*}
	For $w=R_{\tau^n}^{(i)}v$, the form
	$\tau^n \mathfrak a_{t^{(i)},w}(R_{\tau^n}^{(i)}y,R_{\tau^n}^{(i)}z)$
	therefore defines a bounded self-adjoint semi-positive operator on $X$, bounded by $\frac{2C_p^{\mathrm{EIN}}}{\rho^{(i)}}I$. This is the weak interpretation of $\tau^n R_{\tau^n}^{(i)}\mathcal A(t^{(i)},w)R_{\tau^n}^{(i)}$. Apply the square expansions and RK or convex-mixture arguments in the fixed space $X$ separately at each stage.
\end{proof}

Periodic and Neumann cases use their usual form spaces. More general closed semi-positive splitting forms are admissible when their domains support all coefficient and gradient evaluations. The formal temporal order statements above concern fixed spatial discretizations.

\subsection{Systems in a fixed metric}

For a vector unknown $\mathbf u:(0,T] \times \Omega \mapsto \R^r$, consider
\begin{equation}\label{eq:system-model}
	\mathbf u_t-\Div \left(A(t,\mathbf{x},\mathbf u,\Grad \mathbf u)\Grad \mathbf u\right)=\mathbf0.
\end{equation}
Fix $H=H^T\succ0$ and replace the state product by
\begin{equation}\label{eq:system-metric}
	(\mathbf v_h,\mathbf w_h)_H\coloneqq \int_\Omega\mathbf v_h^TH\mathbf w_h\,\d\mathbf{x},\quad
	\norm{\mathbf v_h}_H^2\coloneqq (\mathbf v_h,\mathbf v_h)_H.
\end{equation}
Use the corresponding component metric in the discrete gradient space; for coupled spatial directions use $I_d\otimes H$. The coefficient tensor must be symmetric semi-positive definite in this metric. For a component matrix $A$, this means
\begin{equation}\label{eq:system-coefficient}
	HA=A^TH,\quad H^{\frac{1}{2}}AH^{-\frac{1}{2}}\succeq0.
\end{equation}
For EIN1--3, choose $0\leq C_p^{\mathrm{EIN}}\overline\lambda_A^{n,(i)}(\mathbf{x})\leq s^{n,(i)}(\mathbf{x})$ at every stage, where $\overline\lambda_A^{n,(i)}$ bounds the largest eigenvalue of $H^{\frac{1}{2}}AH^{-\frac{1}{2}}$ evaluated at the filtered vector stage and its discrete gradient at $t^{n,(i)}$. For SSP methods, apply \eqref{eq:stage-substep-bounds} with $a$ replaced by $\overline\lambda_A$ evaluated at the vector branch state. Together with the scalar penalty-weight bounds, this gives \eqref{eq:majorant} in the $H$ metric by the same pointwise integration. For general positive-weight RK methods, use the factor in \eqref{eq:positive-rk-bound}. The stability proofs then apply in $(\cdot,\cdot)_H$ to numerical states $\mathbf u_h^n,\mathbf u_h^{n,(i)}$. In particular, $\norm{\mathbf u_h^{n+1}}_H\leq\norm{\mathbf u_h^n}_H$.

\subsection{Positive-weight quadrature}\label{sec:quadrature}

Choose positive definite quadrature products $(\cdot,\cdot)_Q$ on $V_h$ and $(\cdot,\cdot)_{\Sigma,Q}$ on $\Sigma_h$, using positive weights and injective evaluation maps. For a volume realization,
\begin{equation}\label{eq:quadrature-products}
	(v_h,w_h)_Q \coloneqq \sum_K Q_K(v_h w_h),\quad
	(\mathbf r_h,\mathbf s_h)_{\Sigma,Q}\coloneqq \sum_K Q_K(\mathbf r_h\cdot\mathbf s_h).
\end{equation}
Define $G_{h,Q}$, $\mathcal A_{h,Q}$, and $S_{h,Q}^{n,(i)}$ using these products in the spatial construction, with the same volume and compatible face or patch quadrature for the physical and splitting forms. In particular,
\begin{equation}\label{eq:quadrature-form}
	\begin{aligned}
		(\mathcal A_{h,Q}(t,v_h)z_h,w_h)_Q
		&=\sum_K Q_K\left(a(t,\mathbf{x},v_h,G_{h,Q}v_h)G_{h,Q}z_h\cdot G_{h,Q}w_h\right)\\
		&\quad+c_{h,Q}[\eta_A(t,v_h)](z_h,w_h).
	\end{aligned}
\end{equation}
Define the filters by \eqref{eq:step-indexed-filters} using $S_{h,Q}^{n,(i)}$, and set
\begin{equation}\label{eq:quadrature-energy}
	\begin{aligned}
		w_h&\coloneqq R_{\delta,h,Q}^{n,(i)}v_h,\\
		\mathcal D_{\delta,h,Q}^{n,(i)}(t,v_h)&\coloneqq (\mathcal A_{h,Q}(t,w_h)w_h,w_h)_Q.
	\end{aligned}
\end{equation}
\begin{theorem}[Quadrature realization]\label{thm:quadrature}
	If the coefficient and penalty bounds of \cref{lem:forward,thm:ssp-exact-general,thm:outer} hold at the corresponding quadrature nodes and filtered states, their energy identities and stability estimates hold in $\norm{\cdot}_Q$ with the same stabilization factors. The filter scales remain $\delta=q_{ij}\tau^n$ for SSP branches and $\delta=\tau^n$ for classical RK stages.
\end{theorem}

\begin{proof}
	Positive quadrature weights preserve the coefficient and penalty comparisons, giving \cref{lem:spatial-comparison} in $(\cdot,\cdot)_Q$. The spectral estimates and energy identities then apply in this inner product.
\end{proof}

Exact state mass quadrature gives contraction in the exact $L^2$ norm. Otherwise, mesh-uniform equivalence of $\norm{\cdot}_Q$ and $\norm{\cdot}_{L^2}$ yields a mesh-uniform $L^2$ bound.

\section{Relation to kernel-based resolvent expansions}\label{sec:kernel}

The kernel methods in \cite{10.1007/s10915-020-01152-w,10.4208/cicp.OA-2019-0197} use truncated resolvent expansions with explicit SSP time integration. Christlieb et al. \cite{10.1007/s10915-020-01152-w} write the one-dimensional diffusion term as $\phi(u)_{xx}$, where $\phi'=a$, giving the same Kirchhoff structure as our recent preprint \cite{10.48550/arXiv.2610.34539}. Wang et al. \cite{10.4208/cicp.OA-2019-0197} treat $\Div (A(t,\mathbf{x},u)\Grad u)$, retaining the variable coefficient inside the divergence as in \eqref{eq:model}. These two kernel constructions also differ in filter placement: the first applies one filter to the Kirchhoff potential, whereas the second has a two-sided filter structure in the one-dimensional periodic setting derived below. Thus two-sided filtering is already present in the latter kernel method. Here one-sided and two-sided refer to filter placement relative to the diffusion operator, rather than to the spatial support of the convolution kernel.

For $u_t=\Delta\phi(u)$, let $\mathcal{L}\coloneqq-\Delta$ with compatible boundary conditions and define
\begin{equation}\label{eq:kernel}
	J_\alpha\coloneqq (I+\alpha^{-2}\mathcal{L})^{-1},\quad D_\alpha\coloneqq I-J_\alpha.
\end{equation}
The modified Helmholtz inverse and defect expansion in \cite{10.1007/s10915-020-01152-w} satisfy the exact operator identities
\begin{equation}\label{eq:kernel-identities}
	\alpha^2D_\alpha=J_\alpha \mathcal{L},\quad
	\alpha^2\sum_{j=1}^pD_\alpha^j=\mathcal{L}(I-D_\alpha^p).
\end{equation}
For the corresponding Galerkin operators, let $\phi_h\coloneqq\Pi_{V_h}\phi$. The first partial sum followed by forward Euler is
\begin{equation}\label{eq:kernel-first}
	u_h^{n+1}=u_h^{n}-\tau^n(I+\alpha^{-2}K_h)^{-1}K_h\phi_h(u_h^{n}).
\end{equation}
For $\alpha^{-2}=\tau^n s_0$, this equals the usual single-resolvent EIN1 update in Kirchhoff form with constant splitting coefficient $s_0$. Here $K_h$ is the unit-coefficient discrete Laplacian, and both sides use the same projected potential. In the kernel normalization $\alpha^2=\frac{\beta}{a_{\mathrm{max}}\tau^n}$, the parameter correspondence is $s_0=\frac{a_{\mathrm{max}}}{\beta}$, where $a_{\mathrm{max}}$ is the diffusion bound and $\beta$ is the kernel parameter of \cite{10.1007/s10915-020-01152-w}. This concerns classical EIN1 rather than the two-sided method \eqref{eq:ein1}.

Set $R_\alpha\coloneqq I-D_\alpha^p$. The higher partial sums in \cite{10.1007/s10915-020-01152-w} contain the same filter polynomial as $R_{\tau^n,h}^{(i)}$ in \eqref{eq:resolvent}, after matching scales. Their diffusion action is
\begin{equation}\label{eq:kernel-kirchhoff-filter}
	-\alpha^2\sum_{j=1}^pD_\alpha^j\phi(u)
	=-\mathcal L R_\alpha\phi(u).
\end{equation}
This is one-sided filtering of the potential $\phi(u)$; it does not filter the state before evaluating $\phi$. Its connection with our companion preprint is the Kirchhoff residual and, at first order, the classical EIN1 identity \eqref{eq:kernel-first}. For unit linear diffusion, \eqref{eq:kernel-kirchhoff-filter} gives $-\mathcal L R_\alpha u$, while the two-sided construction gives $-\mathcal L R_\alpha^2u$ at the same filter scale.

For the general divergence structure, consider the one-dimensional periodic operators with exact convolution in Wang et al. \cite{10.4208/cicp.OA-2019-0197}. Let $\mathcal L=-\partial_{xx}$ in \eqref{eq:kernel}, and denote their left and right defects by
\begin{equation}\label{eq:kernel-first-derivative}
	\begin{aligned}
		\mathscr D_L&\coloneqq I-(I+\alpha^{-1}\partial_x)^{-1},\\
		\mathscr D_R&\coloneqq I-(I-\alpha^{-1}\partial_x)^{-1},\\
		\mathscr D_p&\coloneqq\frac{\alpha}{2}\sum_{j=0}^{p-1}
		D_\alpha^j(\mathscr D_L-\mathscr D_R).
	\end{aligned}
\end{equation}
Here their centered defect is $D_\alpha=(\mathscr D_L+\mathscr D_R)/2$. Periodicity gives the exact identities
\begin{equation}\label{eq:kernel-derivative-filter}
	\frac{\alpha}{2}(\mathscr D_L-\mathscr D_R)=J_\alpha\partial_x,
	\quad
	\mathscr D_p=\left(\sum_{j=0}^{p-1}D_\alpha^jJ_\alpha\right)\partial_x
	=R_\alpha\partial_x=\partial_xR_\alpha.
\end{equation}
Their diffusion approximation applies $\mathscr D_p$ to both the solution and the flux. Writing $\mathcal A(t,u)z\coloneqq-\partial_x(a(t,x,u)\partial_xz)$, we therefore obtain
\begin{equation}\label{eq:kernel-two-sided}
	\mathscr D_p\left(a(t,x,u)\mathscr D_pu\right)
	=R_\alpha\partial_x\left(a(t,x,u)\partial_xR_\alpha u\right)
	=-R_\alpha\mathcal A(t,u)R_\alpha u.
\end{equation}
Thus the second kernel paper places $R_\alpha$ on both sides of the frozen physical diffusion operator. Self-adjointness of $R_\alpha$ and skew-adjointness of $\mathscr D_p$ give
\begin{equation}\label{eq:kernel-spatial-dissipation}
	\left(u,\mathscr D_p(a(t,x,u)\mathscr D_pu)\right)
	=-\int_\Omega a(t,x,u)\abs{\partial_xR_\alpha u}^2\,\d x.
\end{equation}

The corresponding residual in the present paper is $-R_{\tau^n,h}^{(i)}\mathcal A_h(t,R_{\tau^n,h}^{(i)}v_h)R_{\tau^n,h}^{(i)}v_h$. Besides sharing the symmetric filter placement in \eqref{eq:kernel-two-sided}, it evaluates the coefficient at the filtered stage, including its filtered discrete gradient when the coefficient depends on $\Grad u$. The kernel construction in \eqref{eq:kernel-two-sided} evaluates $a$ at the unfiltered state $u$. For state-independent coefficients, the two expressions have the same algebraic form when the physical operators and filter scales are matched; for the filters in \eqref{eq:step-indexed-filters} with $\delta=\tau^n$, the scale matching is $\alpha^{-2}K_h=\frac{\tau^n}{2}S_h^{(i)}$. The derivation \eqref{eq:kernel-derivative-filter}--\eqref{eq:kernel-spatial-dissipation} concerns exact one-dimensional periodic convolution.

Both kernel approaches evaluate resolvent inverses through exponential convolutions and fast recurrences, and their constant-coefficient periodic stability analyses use Fourier symbols. For $u_t=a_0u_{xx}$ with $a_0>0$ and $\alpha^2=\beta/(a_0\tau^n)$, the second kernel construction followed by forward Euler has amplification factor $1-\beta x r_p(x)^2$, where $x=\frac{\xi^2}{\alpha^2}$ for Fourier frequency $\xi$. Its contraction condition is $\beta C_p^{\mathrm{EIN}}\leq2$, involving exactly the scalar maximization in \eqref{eq:sharp-ein-constant}. Under our normalization $S=s_0\mathcal L$, matching the filters gives $s_0=2a_0/\beta$, so this condition becomes $0\leq C_p^{\mathrm{EIN}}a_0\leq s_0$, the same threshold as our Euler-branch estimate.

The present construction combines this symmetric filter structure with filtered coefficient evaluation and the stagewise comparison in \cref{lem:forward}. The resulting spectral bound gives fully discrete nonlinear $L^2$ stability through SSP compositions or the RK matrix criterion.

\section{Concluding remarks}\label{sec:conclusion}

We have developed EIN filtered Runge--Kutta methods for nonlinear degenerate diffusion with coefficients depending on space, time, state, and gradient. Applying the same high-order filter to the input and the nonlinear residual gives a symmetric operator product and a non-increasing $L^2$ norm for arbitrary time steps under coefficient and penalty-weight comparisons at the filtered states. The argument allows noncommuting physical and splitting operators, and the splitting operator may vary between stages. The scalar constant $C_p^{\mathrm{EIN}}$ gives the sharp uniform comparison threshold for filtered forward Euler. 

SSP compositions inherit this stability and satisfy exact energy identities accounting for both physical and mixing dissipation. A quadratic matrix criterion extends the construction to explicit RK methods with positive weights, with the additional factor $\gamma_{\mathrm{RK}}$. The examples include SSP-EIN1--3, SSP-EIN(5,4), and classical RK4-EIN. The formal order of the underlying RK method is preserved when order-\(p\) filters are used. Sufficient conditions on the comparison constants for the splitting coefficients are listed in \cref{tab:summary}. Only explicit nonlinear evaluations and symmetric positive definite linear solves are involved. 

The spatial comparison principle applies to CG, BR1, Brezzi, and LDG realizations, and to compatible SBP--SAT discretizations in their discrete mass norm. The same stability argument extends to continuous spatial operators, systems in a fixed metric, and positive-weight quadrature. The comparison with kernel-based resolvent expansions identifies a shared filter polynomial and, for one existing construction, a shared two-sided operator structure. Filtered coefficient evaluation and operator comparison provide the fully discrete nonlinear stability argument for the general diffusion operators considered here.

%\begin{acknowledgements}
%\end{acknowledgements}

\section*{Declarations}

\subsection*{Funding}

No funding is declared for this work.

\subsection*{Conflict of interest}

The authors declare that they have no conflict of interest to disclose.

\subsection*{Data availability}

No datasets were generated or analyzed during the current study.

\bibliographystyle{spmpsci}      % mathematics and physical sciences
\bibliography{refs}   % name your BibTeX data base

\end{document}